\documentclass[journal]{IEEEtran}
\usepackage{amsthm}
\usepackage{outlines}
\usepackage{epsfig}
 \usepackage[english]{babel} 
\usepackage{tabularx,multirow,booktabs,lipsum}
\usepackage{amssymb}
\usepackage{graphics}
\usepackage{graphicx}
\usepackage{subcaption}
\usepackage{epstopdf}
\usepackage{caption}
\usepackage{amsmath, amsthm}
\usepackage{pbox}
\usepackage[numbers]{natbib}
\usepackage[switch]{lineno}
\usepackage{lscape}
\usepackage{longtable}
\usepackage{booktabs}
\usepackage[english]{babel}
\usepackage{color}
\usepackage{etoolbox}
\usepackage{textcomp}
\usepackage{multirow}
\usepackage[normalem]{ulem}

\definecolor{RED}{rgb}{0.7, 0.0, 0.0}
\definecolor{BLUE}{rgb}{0.0, 0.0, 0.7}

\BeforeBeginEnvironment{figure}{\vskip-2ex}
\AfterEndEnvironment{figure}{\vskip-1ex}
\theoremstyle{remark}
\newtheorem{remark}{Remark}

\begin{document}
%
\title{Stochastic and Chance-Constrained Conic  Distribution System Expansion Planning Using Bilinear Benders Decomposition}
%
%
%

\author{Hossein Haghighat,~\IEEEmembership{Member,~IEEE,}
        and~Bo~Zeng,~\IEEEmembership{Member,~IEEE}
}

%
%

\markboth{}%
{}
%



\maketitle

\begin{abstract}
Second order conic programming (SOCP) has been
used to model various applications in power systems, such as operation and expansion planning. In this paper, we present a two-stage stochastic mixed integer SOCP (MISOCP) model for the distribution system expansion planning problem that considers uncertainty and also captures the nonlinear AC power flow. To avoid costly investment plans due to some extreme scenarios, we further present a chance-constrained variant that
could lead to cost-effective solutions. To address the computational challenge, we extend the basic Benders decomposition method and develop a bilinear variant to compute stochastic and chance-constrained MISOCP formulations. A set of numerical experiments is performed to illustrate the performance of our models and computational methods. In particular, results show
that our Benders decomposition algorithms drastically outperform a professional MISOCP solver in handling stochastic scenarios by orders of magnitude.\\
\end{abstract}


%
\IEEEpeerreviewmaketitle

\section*{Nomenclature}
\addcontentsline{toc}{section}{Nomenclature}
\begin{IEEEdescription}[\IEEEusemathlabelsep\IEEEsetlabelwidth{$V_1,V_2,V_3, V_4$}]
\item[$\textbf{\textit{Index and Set}}$]
\item[$\Omega^{E}_{H},\Omega^{C}_{H}$] Set of existing and candidate branches.
\item[$\Omega,  \Omega_{SB}$] Set of network nodes and existing substations.
\item[$\Omega{(i)}$] Set of nodes connected to node \textit{i}.
\item[$S,\Omega_{T}$] Set of scenarios and time blocks representing the target year.
\item[$t,i,j,ij$] Index of time blocks, nodes, and branches.\\
\item[$\textbf{\textit{Parameters}}$]
\item[$c_{f,i}^{sub}, c_{v,i}^{sub}$] Annualized fixed and variable investment costs of substation at node \textit{i}.
\item[$ c_{f,ij}^{fed}, c_{v,ij}^{fed} $] Annualized fixed and variable investment costs of feeder \textit{ij}.
\item[$ c_{o,ij,t,s}^{fed}, c_{o,i,t,s}^{cap} $] Operational costs of feeder \textit{ij} and capacitor at node \textit{i} in scenario \textit{s} and time \textit{t}.
\item[$ c_{f,i}^{cap}, c_{v,i}^{cap}$] Annualized fixed and variable investment costs of capacitor at node \textit{i}.
\item[$ c_{t,s}^{loss}, c_{t,s}^{pn}$] Cost of losses and penalty cost of load curtailment in scenario s and time \textit{t}.
\item[$ g_{ij}, b_{ij},b_{ij}^{sh} $] Series conductance, series susceptance, and shunt susceptance in the \(\pi-\)model of branch \textit{ij}.
\item[$V_{i}^{min}, V_{i}^{max}$] Maximum and Minimum node voltages.
\item[$P_{i}^{smax}, Q_{i}^{smax}$] Active and reactive capacity of substation transformer.
\item[$\alpha_{i}^{s} $]Power factor of substation transformer at node \textit{i}.
\item[$\beta_{i} $] Equals $\tan{(cos^{-1}{\phi_i})}$ with $\phi_i$ being the power factor of the load at node \textit{i}.
\item[$\pi_{s}$] Probability of scenario \textit{s}.
\item[$d_{t}^{h}$] Number of hours comprising time block \textit{t}.
\item[$Q_{c,i}^{max} $] Capacity of capacitor invested at node \textit{i}.
\item[$  l_{ij}^{fed},I_{ij}^{max}$] Feeder length and rated current.
\item[$ P_{i,t,s}^{d},Q_{i,t,s}^{d}$] Active and reactive demand at node \textit{i}.\\
\item[$\textbf{\textit{Variables}}$]
\item[$v_{i}^{sub},v_{i}^{cap} $] Binary variable which equals 1 if investment is made for substation/capacitor and zero otherwise.
\item[$ k_{ij}$] Binary variable which equals 1 if investment is made for any (new/existent) branch \textit{ij} and is zero otherwise.
\item[$ f_{ij}$] Binary variable associated with the existent branch \textit{ij} which equals 1 if utilized and is zero otherwise.
\item[$z_{ij} $]Binary variable which is equal to 1 if \textit{j} is the parent node of \textit{i}.
\item[$y_{ij} $] Connection status of branch \textit{ij}: 1 connected; 0 disconnected.
\item[$I_{ij}^{fed} $] Current flow of branch \textit{ij}.
\item[$ S_{i}^{sub} $] Investment in substation at node \textit{i}.
\item[$P_{ij,t,s},Q_{ij,t,s} $]Real and reactive flow on branch \textit{ij}.
\item[$ P_{Ii,t,s},Q_{Ii,t,s} $] Real and reactive injections at node \textit{i}.
\item[$ P_{i,t,s}^{sb}, Q_{i,t,s}^{sb}$] Real and reactive injections at substation \textit{i}.
\item[$Q_{i}^{cap} $] Reactive injection of capacitor at node \textit{i}.
\item[$r_{i,t,s} $] Real power curtailed at node \textit{i}.
\item[$u_{i,t,s}^{ij}$, $u_{i,t,s}$] Voltage of node \textit{i} associated with branch \textit{ij}.
\item[$C_{inv}, C_{opr,s} $] Investment cost and operation cost of scenario \textit{s}.
\item[$M$] A sufficiently large positive number.
\end{IEEEdescription}

\section{Introduction}
\IEEEPARstart{T}he aim of distribution system expansion planning is to derive a cost-effective investment plan on system configuration. The problem is often formulated as an optimization model that minimizes a desired objective such as installation costs of  new facilities (e.g., feeders, substation, etc.), the upgrade costs of existing equipments, operational costs, while respecting various specifications or operating limits. In the context of mathematical programming, mixed-integer linear programs, e.g., \cite{Vaziri} and \cite{Haffner}, are primarily constructed, given that they are flexible in modeling and computationally well supported by commercial solvers. Nevertheless, power flow equations are inherently nonlinear and simple linear approximations may lead to solutions of bad qualities. Hence, to capture critical factors that are beyond the capacity of linear models, nonlinear formulations are employed. Among them, it is worth mentioning the convex second order conic programming (SOCP) model in \cite{jabr2008optimal} to describe nonlinear AC power flow equations, which actually is an exact model for a radial distribution system. Note that SOCP model is a convex optimization formulation that can be efficiently computed by professional solvers. Moreover, because discrete investment decisions are generally involved in system planning, integer variables are introduced and that formulation is extended to a mixed-integer SOCP (MISOCP) model \cite{jabr2013polyhedral}. Recent applications of conic optimization in power system studies include \cite{Taylor},\cite{Lee} on network reconfiguration for loss reduction, and \citep{Ding_SOCP} on reactive power optimization. We observe that heuristic techniques,e.g.,  \citep{Heurist_1,Heurist_2} have also been applied to the problem of distribution expansion planning. However,  optimality of their solutions is not guaranteed.\par
One fundamental challenge of a planning problem is the uncertainty, as it is for the future situation and perfect forecasts on demands or cost parameters are impossible to obtain. To include the uncertainty in decision making, the most popular approach is to build a two-stage stochastic program model \cite{birge2011introduction}. Under this modeling scheme, the first-stage decisions  are those made before the realization of random factors, e.g., investment decisions; the second stage decisions are those
made after observing the actual realization of those factors, e.g., actual operations.
Using a practical strategy to represent the randomness by a finite number of scenarios (and their realization probabilities), one set of second stage decision variables are introduced for every single scenario, which therefore leads to a large-scale mathematical program. Very often, the decision maker notes that some random situations might be rare but very costly to handle. With that observation,  he may just want to adopt a cost-effective plan that protects himself from most random scenarios, e.g., 95\% of all possible situations, and ignore the remaining ones. To meet this modeling demand, the original  stochastic program could be extended to the chance-constrained program \citep{ahmed2008solving}. Specifically, a binary variable is introduced to associate with each scenario, which enables the decision makers to select his concerned scenarios.
Through using scenario-based stochastic or chance constrained modeling schemes, many practical decision making problems in power, transportation or healthcare systems, which are subject to serious random factors, have been successfully addressed \citep{birge2011introduction,ahmed2008solving}. Typical power systems applications include unit commitment problems \citep{Takriti_1996,wang2012chance,Zheng_2015}; generation and transmission  planning and capacity expansion \citep{Gil_2014,parpas2014stochastic}; storage siting and sizing \citep{baker2014optimal,kuznia2013stochastic}; and grid vulnerability analysis \citep{cormican1995computational,janjarassuk2008reformulation}. Note that the majority of the developed models are (stochastic) mixed integer linear programs. To handle the large-scale structures from stochastic scenarios, various Benders decomposition algorithms have been customized and developed \citep{benders1962partitioning,van1969shaped,shahidehopour2005benders}, which could drastically reduce computational times and render those models useful for practitioners.\par
According to the recent survey on \citep{ganguly2013recent}, however, published works incorporating uncertainty into the distribution system planning problem are few and far between.  Except a few heuristic methods \citep{Heurist_1,Heurist_2,Heurist_3} to consider the  uncertainty issue, there is no analytical method to systematically support the practice of distribution system planning under random environments. One essential challenge is the large number of nonlinear conic AC power flow equations associated with stochastic scenarios. Such large-scale MISOCP formulation is generally beyond the solution capacity of commercial mixed integer conic solvers. Indeed, except a theoretical analysis in \citep{wei2015generalized}, we have not observed any efficient decomposition method for large-scale MISOCP problems. In \citep{wei2015generalized},  authors presents a generalized Benders decomposition method that theoretically converges to an optimal solution of one type of mixed integer conic problem. Neither its performance over practical instances nor its extension to deal with stochastic conic programs is reported or analyzed.
 To  address such a situation in the research of distribution system planning, in this paper, we first propose a stochastic two-stage MISOCP formulation and its chance-constrained variant, where the first stage is the investment stage that determines the appropriate conductor types and feeder construction routs together with substation reinforcement and capacitor installation. The second stage is the operating stages where operations decisions are derived to minimize costs of involved losses and maintenance. The uncertainty of loads and energy price is modeled through a set of scenarios with their realization probabilities. Because both the stochastic and chance-constrained models are challenging to commercially available solvers, we design and implement fast Benders decomposition methods that can efficiently handle a large number of  stochastic scenarios containing conic flow equations. Our major contributions are listed as follows:\par
 $(i)$ On the modeling aspect, we develop the first two-stage stochastic mixed integer conic model for the distribution system expansion planning that considers uncertainty and accurately captures the nonlinear AC power flow. We also present its chance-constrained variant to avoid the costly investment plans due to some extreme scenarios.\par
$(ii)$ On the algorithm development aspect, we develop fast Bender decomposition algorithms that drastically outperforms a professional mixed integer conic solver in handling stochastic scenarios  by orders of magnitude. To the best of our knowledge, they are the first Benders decomposition methods that can practically deal with stochastic and chance-constrained mixed integer conic programs. Indeed, our customized Benders decomposition methods are general approaches that can be widely applied to address many other applications.\par
$(iii)$ On the computation aspect, we perform a set of preliminary experiments. Results show that  investment plans derived from stochastic and chance-constrained models
present better performances, under the stochastic environments, over those derived from the deterministic counterparts. Results also highlight that our developed algorithms drastically outperforms professional conic programming solvers and can derive optimal solutions in a reasonable time.\par

The paper is organized as follows. In Section II, we develop the mathematical formulation of the two-stage stochastic distribution system expansion planning problem and its chance-constrained variant. In Section III, we propose our Benders decomposition methods, including the master and subproblem formulations, and the algorithmic operations. In Section IV, computational experiments are presented and analyzed. Finally, relevant conclusions are drawn in Section V.
\section{Stochastic and Chance-constrained Distribution System Expansion Planning Models}
In this section, we present the stochastic MISOCP distribution system expansion model with detailed descriptions, followed by its chance-constrained variant.

\subsection{Stochastic distribution system expansion planning model}
\label{sect_SP_full}
We consider the distribution system expansion planning under a two-stage scheme, where the first stage determines the expansion decisions on equipment and the second stage models system operations with the system upgrades. Given that operational decisions are made according to real time load and electricity price, which are random, we adopt a set of discrete scenarios (and their realization probabilities) to represent their uncertainties and define the second stage  recourse problem for every scenario. Because of the randomness in loads, we allow possible load curtailment in those scenarios,  which will then be penalized in the cost function of each recourse problem. As a result, we have a two-stage stochastic programming formulation for distribution system expansion planning problem, as in the following.
\begin{subequations}
\begin{align}
\min &  \ \ {C_{inv}} + \sum\limits_{s \in S} {{\pi _s}} {C_{opr,s}} \label{eq_obj}\\
\mbox{s.t.} \ {C_{inv}} =& \sum\limits_{i \in {\Omega _{SB}}} {\left( {c_{f,i}^{sub}v_i^{sub} + c_{v,i}^{sub}S_i^{sub}} \right)}  \\
+& \sum\limits_{ij \in \Omega _H^C \cup \Omega _H^E} c_{f,ij}^{fed}l_{ij}^{fed}{k_{ij}} \nonumber\\
+ & \sum\limits_{i \in \Omega \backslash {\Omega _{SB}}} {\left( {c_{f,i}^{cap}v_i^{cap} + c_{v,i}^{cap}Q_i^{cap}} \right)} \nonumber \\
 \ \ {C_{opr,s}} =& \sum\limits_{t \in {\Omega _T}} {d_t^h}\Big(\sum\limits_{i \in \Omega} {c_{t,s}^{pn} r_{i,t,s}+c_{t,s}^{loss}{P_{Ii,t,s}}} \nonumber \\
 + & \sum\limits_{ij \in \Omega _H^C \cup \Omega _H^E} {c_{o,ij,t,s}^{fed}{y_{ij}}}\Big)  \label{eq_sce_opt}\\
 P_{Ii,t,s} & = P_{i,t,s}^d - r_{i,t,s} + \sum_{j \in \Omega (i)} {{P_{ij,t,s}}}  \,\,\forall i \in \Omega ,t,s \label{eq_a2}\\
P_{i,t,s}^{sb} &= \sum\limits_j {{P_{ij,t,s}}} \qquad\forall i \in {\Omega _{SB}},t,s \label{eq_a3}
\end{align}
\begin{align}
Q_{Ii,t,s}  &= Q_{i,t,s}^d - \beta_i r_{i,t,s} + \sum\limits_{j \in \Omega (i)} Q_{ij,t,s} \, - Q_i^{cap}  \nonumber \\& \qquad \forall i \in \Omega ,t,s
\label{eq_a4}\\
Q_{i,t,s}^{sb} &= \sum\limits_j {{Q_{ij,t,s}}} \, \qquad\forall i \in {\Omega _{SB}},t,s \label{eq_a5}
\end{align}
\begin{align}
\noindent \sqrt 2 ( g_{ij}^2 & + (b_{ij}^{} + 0.5b_{ij}^{sh})^2) u_{i,t,s}^{ij} + \sqrt 2 \left( {g_{ij}^2 + b_{ij}^2} \right)u_{j,t,s}^{ij} \nonumber\\
- 2( g_{ij}^2 & + b_{ij}^2 + 0.5b_{ij}^{sh}b_{ij}) {R_{ij,t,s}} + g_{ij}b_{ij}^{sh}{L_{ij,t,s}} \nonumber\\
\le & (I_{ij}^{fed})^2 \qquad \forall ij \in \Omega _H^E \cup \Omega _H^C,t,s  \label{eq_a6}\\
0 \le u_{i,t,s}^{ij}& \le \frac{{{{(V_i^{\max })}^2}}}{{\sqrt 2 }}{y_{ij}}  \qquad\forall i \in \Omega ,t,s  \label{eq_a7}\\
0 \le u_{j,t,s}^{ij}& \le \frac{{{{(V_j^{\max })}^2}}}{{\sqrt 2 }}{y_{ij}}  \qquad\forall j \in \Omega ,t,s  \label{eq_a8}\\
0 \le {u_{i,t,s}}&\!-\!u_{i,t,s}^{ij} \le \frac{{{{(V_i^{\max })}^2}}}{{\sqrt 2 }}(1 - {y_{ij}})  \quad\forall i \in \Omega ,t,s  \label{eq_a9}\\
0 \le{u_{j,t,s}} &\!-\!u_{j,t,s}^{ij} \le \frac{{{{(V_j^{\max })}^2}}}{{\sqrt 2 }}(1 - {y_{ij}}) \quad\forall j \in \Omega ,t,s \label{eq_a10}\\
0 \le{R_{ij,t,s}}& = {R_{ji,t,s}}  \qquad\forall ij \in \Omega _H^E \cup \Omega _H^C,t,s \label{eq_a11}\\
\frac{{{{(V_i^{\min })}^2}}}{{\sqrt 2 }} &\le{u_{i,t,s}} \le \frac{{{{(V_i^{\max })}^2}}}{{\sqrt 2 }} \qquad\forall i \in \Omega ,t,s  \label{eq_a12}\\
{L_{ji,t,s}} =  - &{L_{ji,t,s}}  \qquad\forall ij \in \Omega _H^E \cup \Omega _H^C,t,s \label{eq_a13}\\
P_{i,t,s}^{sb} \le & P_i^{{\rm{smax}}} + \alpha_i^sS_i^{sub} \quad\forall i \in {\Omega_{SB}},t,s  \label{eq_a14}\\
Q_{i,t,s}^{sb} \le & Q_i^{{\rm{smax}}} +\sin(\cos^{-1}(\alpha_i^s))S_i^{sub} \ \forall i \in {\Omega_{SB}},t,s \label{eq_a15}\\
0 \le S_i^{sub} & \le Mv_i^{sub}   \qquad\forall i \in {\Omega_{SB}}   \label{eq_a16}\\
0 \le Q_i^{cap} & \le v_i^{cap}Q_{c,i}^{\max }   \qquad\forall i \in \Omega \backslash {\Omega _{SB}}                \label{eq_a17}\\
{k_{ij}} + {f_{ij}} & = 1      \qquad\forall ij \in \Omega _H^E \cap \Omega _H^C            \label{eq_a18}\\
{y_{ij}} \le {f_{ij}} &        \qquad\forall ij \in \Omega _H^E - \Omega _H^E \cap \Omega _H^C         \label{eq_a19}\\
{y_{ij}} \le {k_{ij}} &       \qquad\forall ij \in \Omega _H^C - \Omega _H^E \cap \Omega _H^C      \label{eq_a20}\\
{z_{ij}} = 0 &         \qquad\forall ij \in \Omega _H^E \cup \Omega _H^C,i \in {\Omega _{SB}}         \label{eq_a21}\\
{z_{ij}} + {z_{ji}} & = {y_{ij}}      \qquad\forall ij \in \Omega _H^E \cup \Omega _H^C            \label{eq_a22}\\
\sum_{j: ij \in \Omega _H^E \cup \Omega _H^C}&{{z_{ij}}}   = 1   \,\,\qquad\forall i \in \Omega \backslash {\Omega _{SB}} \label{eq_a23}\\
(I_{ij}^{fed})^2 = &   (I_{ij}^{\max })^2(1 - {k_{ij}})      \qquad\forall ij \in \Omega _H^E         \label{eq_a24}\\
(I_{ij}^{fed})^2 = & (I_{ij}^{\max })^2 k_{ij}         \qquad\forall ij \in \Omega _H^C     \tag{1aa}    \label{eq_a25}
\end{align}
\begin{align}
2u_{i,t,s}^{ij}u_{j,t,s}^{ij} &\ge R^2_{ij,t,s} + L^2_{ij,t,s}      \ \forall ij \in \Omega _H^E \cup \Omega _H^C,t,s    \tag{1ab}  \label{eq_a26}
\end{align}
\begin{align}
 P_{ij,t,s}  = \sqrt 2 & g_{ij}u_{i,t,s}^{ij} - {g_{ij}}{R_{ij,t,s}} - b_{ij}L_{ij,t,s}  \nonumber\\
 \forall ij \in \Omega_H^E & \cup \Omega_H^C,t,s    \tag{1ac}  \label{eq_a27}\\
Q_{ij,t,s} =  - \sqrt 2 &   (b_{ij} + b_{ij}^{sh}/2)u_{i,t,s}^{ij} + b_{ij}R_{ij,t,s} - {g_{ij}}{L_{ij,t,s}} \nonumber \\
\forall ij \in \Omega_H^E & \cup \Omega_H^C,t,s \tag{1ad}  \label{eq_a28}\\
 v_i^{sub},v_i^{cap}, & {y_{ij}}, \,\,{k_{ij}},{f_{ij}},{z_{ij}}\,\,\,{\rm{binary}},\,\,r_{i,t,s} \ge 0    \tag{1ae}        \label{eq_a29}
\end{align}
\end{subequations}

In \eqref{eq_obj}-\eqref{eq_sce_opt}, the objective is to derive an optimal expansion plan to minimize the investment cost, including both fixed costs and variable costs of new facilities in
 substations, feeders and capacitors, in the first stage and the expected operational cost obtained from every scenario $s\in S$ in the second stage, including cost of losses and  maintenance cost of installed  branches. Constraints in \eqref{eq_a2}-\eqref{eq_a13} are the second order conic description of AC power flow model on a radial distribution system. The capacity limit of the substation is represented by \eqref{eq_a14}-\eqref{eq_a15}. Constraints in \eqref{eq_a16}-\eqref{eq_a17} are investment decisions on substation and capacitors banks.
 The logical constraints \eqref{eq_a18}-\eqref{eq_a20} state that if any investment is made (i.e., $k_{ij}$=1) for an existing branch (by replacing it with a conductor of a higher capacity), the old  conductor must be disconnected (i.e. $y_{ij}$=0) from the network  considering the radiality requirement on the network topology.  Constraints in \eqref{eq_a21}-\eqref{eq_a23} enforce the radiality of the expanded network. Constraints in \eqref{eq_a24}-\eqref{eq_a25} provide capacities of new and existing branches. Variable definitions are presented in \eqref{eq_a26}-\eqref{eq_a29}.

Existing commercial MISOCP solvers can compute small-scale instances. Nevertheless, as shown in Section \ref{sect_numerical}, with the size of distribution network or the number of stochastic scenarios increased, it is necessary to develop advanced algorithms. To this end, a customized Benders decomposition for such MISOCP is developed in Section \ref{sect_benders}.

\subsection{Chance-constrained SOCP Model}
 To introduce our chance-constrained model without repetitive information, we
 represent stochastic programming \textbf{StoP} model in the following compact form.  We use
  $x$ to denote the first stage (investment) variables and
  $y_s$ to denote the second stage (operation) variables in scenario $s$, which is associated with parameters $(g_s,E_s,d_s,B_s, l_s, H_s,h_s)$. Note that single variable restrictions are in \eqref{eq_a29}.
\begin{subequations}
\begin{align}
 {\bf{StoP}}{\rm{: }} \ \min cx &+ \sum_{s\in S}\pi_sg_sy_s    \label{eq_b1}\\
s.t. \ & Fx \le f        \label{eq_b2}\\
 & E_sy_s = d_s   \qquad \forall s  \label{eq_b3}\\
 & Ax + B_sy_s \ge l_s  \qquad \forall s  \label{eq_b4}\\
 & \left\| {H_sy_s} \right\| \le h_sy_s   \qquad \forall s \label{eq_b5}
\end{align}
\end{subequations}\vspace{-2pt}
  Compared to the complete
  formulation in Section \ref{sect_SP_full}, constraints in \eqref{eq_b2} are those for investment decision variables,  constraints in \eqref{eq_b3} include power balance equations, constraints in \eqref{eq_b4} link investment and operation decisions, and constraints in \eqref{eq_b5} represent the second order cone equations. Note that the second term in the objective function \eqref{eq_b1} computes the expected cost over all scenarios.  Next, we provide its chance-constrained variant using the conventional Big-M method. \vspace{-2pt}
    \begin{subequations}
\begin{align}
 {\bf{CC\!-\!bigM}}: \ & \ \min \ cx + G(y_1,w_1,\dots,y_{|S|},w_{|S|})    \label{eq_d1}\\
    s.t. \qquad  & Fx \le f         \label{eq_d2}\\
 & E_sy_s-d_s +Mw_s \ge 0 \qquad \forall s  \label{eq_d3}\\
  & {d_s} - {E_s}{y_s} + M{w_s} \ge 0 \qquad \forall s  \label{eq_d4} \\
  & Ax + {B_s}{y_s} + M{w_s} \ge {l_s} \qquad \forall s \label{eq_d5}\\
 & \left\| {{H_s}{y_s}} \right\| - M{w_s} \le {h_s}{y_s} \qquad \forall s \label{eq_d6}\\
 & \sum\limits_{s = 1}^{|S|} {{\pi _s}{w_s} \le \varepsilon }         & \label{eq_d7}\\
 & {w_s} \in \{ 0,1\} \ \forall s\in S.
\end{align}
\end{subequations}
Clearly, if $w_s=1$, all constraints in scenario $s$ can be ignored
due to Big-M. Hence, binary variable $w_s$ can be used to reflect the inclusion of scenario $s$
in computing an optimal solution.
According to \eqref{eq_d7}, we
seek for a solution that performs well in $(1-\varepsilon)\times 100\%$ of all random situations. Given that some scenarios are ignored in the solution evaluation, we introduce $G$, a function of $y_s$'s and $w_s$'s,  in the objective function \eqref{eq_d1},  to capture the cost contribution of the second stage decisions.

Although the aforementioned $\mathbf{CC - bigM}$ formulation can be treated as
regular mixed integer program if $G$ is defined appropriately, it is noted in \citep{Zeng_An_Kuznia_2014}
that its computational burden is very heavy. For the case that we only care about costs incurred in the concerned
scenarios, we next provide a  bilinear reformulation that
typically has a better computational performance than the Big-M based formulation \citep{Zeng_An_Kuznia_2014}. Moreover, such bilinear format allows us to generalize Benders decomposition method to further improve our solution capacity.
\begin{subequations}
 \begin{align}
{\bf{CC\!-\!BL}}: & \min \,\,\,cx + \rho \sum\nolimits_{s = 1}^{|S|} \pi_s\eta{_s}            \label{eq_e1}\\
      s.t.\,\,      & Fx \le f                       \label{eq_e2}\\
     & \eta_s = g_sy_s( 1 - {w_s} )    \qquad  \forall s   \label{eq_BLe3}\\
     & \left( {{E_s}{y_s} - {d_s}} \right)\left( {1 - {w_s}} \right) = 0 \ \qquad   \forall s  \label{eq_e4}\\ 
     & \left( {Ax + {B_s}{y_s} - {l_s}} \right)\left( {1 - {w_s}} \right) \ge 0 \  \qquad   \forall s  \label{eq_e5}\\ 
     & \left( {\left\| {{H_s}{y_s}} \right\| - {h_s}{y_s}} \right)\left( {1 - {w_s}} \right) \le 0   \  \qquad \forall s   \label{eq_e6}\\ 
     & \sum\limits_{s = 1}^S {{\pi _s}{w_s} \le \varepsilon }                      \label{eq_e7}\\
     & {w_s} \in \{ 0,1\}  \qquad         \forall s \label{eq_e8}
\end{align}
\end{subequations}

 Note from \eqref{eq_BLe3} that if $w_s$ is set to one, i.e., scenario $s$ is ignored, it dose not contribute to the total cost. So, by assigning $w_s$ to one or zero, the impact of scenario $s$, including the cost contribution of recourse decisions and feasibility requirements from recourse constraints, will be explicitly removed from or imposed in the whole formulation. Parameter $\rho$ is introduced to reflect our attitude towards the recourse cost from the concerned scenarios.

\begin{remark}If $\varepsilon=0$, we have $w_s=0$ for all $s\in S$, which reduces $\mathbf{CC\!-\! bigM}$ or $\mathbf{CC\!-\!BL}$ formulations to corresponding \textbf{StoP} model (with the second stage cost weighted by $\rho$).
\end{remark}

\section{Customized Bilinear Benders Decomposition for Stochastic  MISOCP Models}
\label{sect_benders}
Benders decomposition, which is a master-subproblem structured method, is probably the most popular approach to compute scenario-based stochastic mixed integer linear programs. Nevertheless, there is little study on extending this classical method to compute stochastic MISOCP. Actually, we observe in a related study \citep{Ding_SOCP}  that robust MISOCP can be solved by a similar master-subproblem structured decomposition method. With this observation, we extend and customize the basic Benders decomposition method to compute stochastic MISOCP formulation and its chance constrained variant. As stochastic MISOCP formulation can be treated as a special one of its chance constrained variant, we describe our Benders decomposition method in the context of chanced constrained model, particularly the bilinear $\mathbf{CC-BL}$ form,  to simplify our exposition. We next present the subproblem and the master problem, and the solution algorithm.

\subsection{Subproblem and master problem}
Note that we allow demand curtailment in the second stage recourse problem, which is penalized in the objective function. Hence, the strict feasibility of the recourse problem is ensured, which actually guarantees the strong duality of this SOCP problem.
Next, we define the subproblem that is constructed by taking the duality of the second stage recourse problem in scenario $s$ in the $i^\textrm{th}$ iteration for given first stage $\hat x^i$.
\begin{subequations}
\begin{align}
\mathbf{SP}_s: J_s = \max\,\,\,& {\lambda_s}\left(l_s - A \hat x^i \right) + {d_s}{\theta_s}   \\
   s.t.\,\, &   {E_s}{\lambda _s} + {B_s}{\theta _s} + {H_s}{\sigma _s} + {\mu _s}{h_s} = g{ _s}   \\
   & \parallel\sigma_s\parallel \le {\mu _s}        \\
   & \theta_s, \mu_s \ge 0,\,\,\,\, \lambda_s, \sigma_s \,\, \rm{free}
\end{align}
\end{subequations}
Note that the second stage recourse problem and $\mathbf{SP}_s$ always have a finite optimal value.  Hence, we can derive an optimal solution to $\mathbf{SP}_s$ that is an extreme point of its feasible region. Let $(\hat \theta^{i}_s, \hat \mu^i_s, \hat \lambda^i_s, \hat \sigma^i_s)$ denote that optimal solution.

Next, we define the master problem for the $i^\textrm{th}$ iteration. Note that the conventional Benders decomposition method simply generates Benders cuts that are linear functions of $(\hat \theta^{i}_s, \hat \mu^i_s, \hat \lambda^i_s, \hat \sigma^i_s)$. Different from that, we follow  $\mathbf{CC-BL}$ to modulate our Benders cuts by $w_s$ in a bilinear form.
\begin{subequations}
\begin{align}
{\bf{MP\!-\!BL}}: & \min \,\,\,cx + \rho \sum_{s\in S} \pi_s\eta{_s} \label{eq_e1} \\
s.t. \,\,      & Fx \le f                       \label{eq_e2}\\
  & \Big(\hat \lambda^k_{s}(l_{s}-Ax)+d_{s}\hat \theta^k_{s} \Big)(1-w_{s}) \leq \eta_{s}, \nonumber\\
 & 1\leq k<i, \ s\in S  \label{eq_e3}\\
    & \sum_{s \in S} {\pi _s}{w_s} \leq \varepsilon   \label{eq_e4} \\
  & w_s \in \{ 0,1\}  \qquad    s\in S\label{eq_e5}
\end{align}
\end{subequations}
We mention that by enumerating all extreme points of $\mathbf{SP}_s$ (which could be infinite due to the conic structure), $\mathbf{MP-BL}$ is the exact Benders reformulation of $\mathbf{CC-BL}$. Hence, a formulation of $\mathbf{MP-BL}$ defined on a subset of extreme points is a relaxation to $\mathbf{CC-BL}$, whose optimal value provides a lower bound. Note that the bilinear terms in  \eqref{eq_e3} can be linearized using McCormick linearization technique, which converts $\mathbf{MP-BL}$ into an MIP.

\subsection{Solution algorithm}
Obviously, any feasible solution to $\mathbf{CC-BL}$ provides an upper bound. So, iteratively including optimality cuts \eqref{eq_e3} could help us generate better lower and upper bounds. Let $LB$ and $UB$ be the current lower and upper bounds, and $e$ be the optimality tolerance. \\
\hrule\vspace{2mm}
\textbf{Bilinear Benders Decomposition for $\mathbf{CC-BL}$}
\hrule\vspace{2mm}
\begin{outline}[enumerate]
\1 \textbf{Initialization}: Set $LB = -\infty$, $UB= \infty$ and the iteration counter $i=0$.
\1 \textbf{Iterative steps}:
    \2 Compute the master problem $\mathbf{MP-BL}$. Derive its optimal value, $V_i$, and an optimal solution $(\hat x^i, \hat w^i)$. Update lower bound $LB= V_i$.
    \2 For  every $s\in S$, compute subproblem $\mathbf{SP_s}$, derive its optimal value $\hat J_s$ and generate a Benders cut  \eqref{eq_e3} based on its optimal solution.
    \2 Compute $\bar J = \min \{\displaystyle\sum_{s\in S} \pi_s\hat J_s(1-w_s): \displaystyle\sum_{s\in S}\pi_sw_s\leq \epsilon\}$ and update upper bound $UB= \min\{UB, c\hat x^i+\bar J\}$.
\1 \textbf{Stopping criteria}: if $|\frac{UB-LB}{LB}| \le e$, terminate with a solution associated with $UB$. Otherwise, $i=i+1$ and go back to step 2.

\vspace{2mm}
\hrule\vspace{2mm}
\end{outline}
When $\epsilon=0$, our chance-constrained model reduces to the stochastic programming formulation. Because all scenarios must be considered, we can eliminate $w_s$ variables from our bilinear Benders decomposition.

\begin{remark}
Subproblem $\mathbf{SP}_s$ be converted into a linear program using a polyhedral relaxation of the convex feasible region suggested in \citep{Ben-Tal}. We observe, however, that if doing so, the size and runtime of the subproblem increase substantially. The strength of our procedure is that it can attain the exact solution in a fast way than using the polyhedral approximated subproblem, which actually yields an approximate solution and increases the problem dimensionality.
\end{remark}

\section{Numerical Results}
\label{sect_numerical}
In this section a simple example and a test system are employed to demonstrate the application of the proposed methods. They are implemented in AMPL with the optimality tolerance $e$ being $0.1\%$ for Benders decomposition, and with the most popular professional solver CPLEX \citep{Fourer} under default settings. Our experiments are performed on a Windows-based PC with two 3 GHz processors and 4 GB RAM, with time limit set to 3,600 seconds.
\subsection{Example}
For illustration purpose, we applied proposed the methods to the example system shown in Fig.~\ref{fig_1}. This is an originally meshed network to be expanded for target year 5. Load data is provided in Table~\ref{table_I}. Fixed and variable investment costs of substation were assumed 200000 US\$, and 50000 US\$/MW. Feeders investment and maintenance costs were assumed 150000 US\$/km and 450 US\$, respectively \citep{Bhattacharya}. Capacitors fixed and variable investment costs were 3000 US\$, 450 US\$/kVAr  \citep{Chen}. Annualized costs were computed assuming a life span of 15 years with an interest rate of 10\%. The cost of losses was assumed to be ten times the prices given in Table~\ref{table_I}.
We considered seven candidate branches, involving replacements for existent branches (with higher capacities) and new candidate branches 5-4 and 3-5. Results of two formulations denoted as \textit{chance-unconstrained} and \textit{chance-constrained} are presented next.

\subsubsection{Chance-unconstrained case}
We solved three cases involving a deterministic case (denoted as case 0) and two stochastic cases (denoted as case 1 and case 2) where load and price are subject to uncertainty. To describe uncertainty, we created ten equiprobable scenarios (for each case) in which price and load are scaled according to the data provided in Table~\ref{table_II}.
The expansion results are shown in Fig.~\ref{fig_1}. In case 0 no investment is made and a radial topology is achieved by removing branch 3--4. Results in the stochastic case 1 indicate that two existent branches 1--2 and 2--5 are replaced with branches of higher capacity and the substation is expanded by 0.110 pu. In case 2, no investment is required and the existent branch 2--4 is removed to enforce radiality. We observe that different expansion plans are achieved as a result of incorporating uncertainty in the problem.

\subsubsection{Chance-constrained case}
\begin{figure}[t]
\centering
\includegraphics[width=8cm, height=8cm]{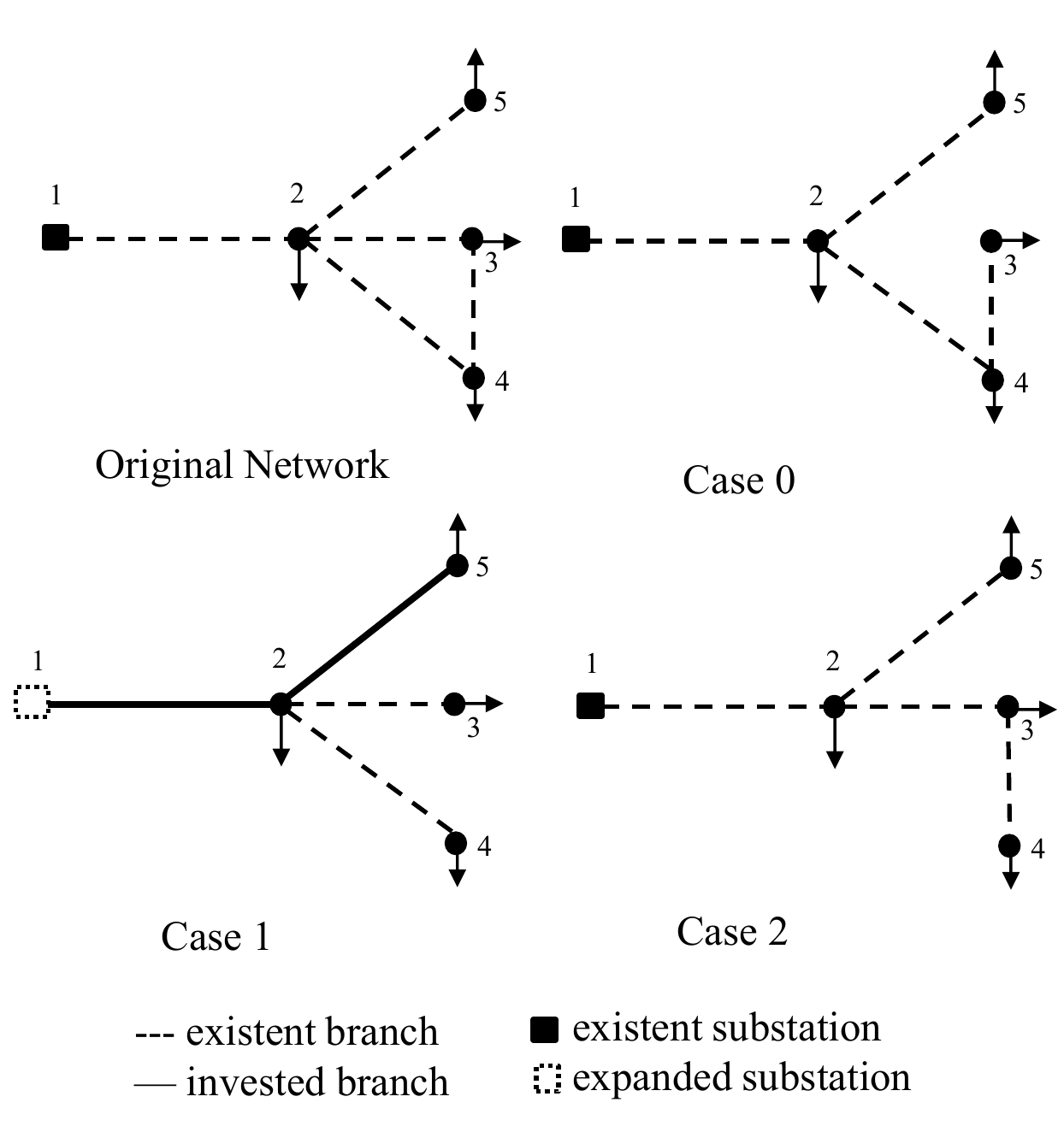}
\caption{Original and expanded networks}
\label{fig_1}
\end{figure}

In this experiment the impact of chance constraint is explored. The \textit{chance-constrained} model was solved for scenarios described in case 1 of Table~\ref{table_II}. Results for chance  levels (namely, 1-$\varepsilon$) 100, 90, 80, and 65 percent are presented in Table~\ref{table_III}. For this small system, they actually lead to a configuration same as that of case 1 in Fig.~\ref{fig_1}. Observe that the most expensive scenarios, namely 4, 3, and 6 (see Table~\ref{table_II}), are dropped when the chance level decreases. Hence, by adjusting the chance level,  the decision maker can have a
balance between the cost of the expansion plan (which results in a high capacitated system) and the desired level of security against risks.

\subsubsection{Computational performance}
In Table~\ref{table_IV}, we compare the computational performance of following different methods. If an optimal solution cannot be derived due to the time limit, the instance is labeled with ``T'' and the corresponding optimality gap, if available, is reported.
 \begin{enumerate}
  \item \textit{Non-decomposed-SOCP}: Simply using MISOCP solver to  compute;
  \item \textit{Benders-linear StoP}: Using Benders decomposition method, with the subproblem being the linearized approximation to the second stage SOCP recourse problem through the method proposed by \citep{Ben-Tal};
  \item \textit{Benders-SOCP StoP}: Using standard Benders decomposition for stochastic programming model;
  \item \textit{Benders-SOCP CC}: Using the bilinear Benders decomposition for chance-constrained model.
\end{enumerate}
Based on Table~\ref{table_IV}, it is obvious that the state-of-the-art CPLEX solver
does not have the capacity to handle practical instances. For this 5-node system, which is definitely small-scale, CPLEX fails to produce any feasible solution for instances with 50 or 100 scenarios within 3,600 seconds. On the contrary, Benders decomposition methods show a drastically improved solution capacity, which reduce the computational time by orders of magnitude. Specifically, Benders methods, which directly call solver to compute subproblems, can compute all instances within a short time. Although the chance constrained formulations could be much more challenging to compute than their stochastic programming counterparts, we observe that our bilinear Benders method demonstrates a strong power to handle chance-constrained second order conic programs,
which is comparable to the classical Benders method for pure stochastic programs. Indeed, within a couple of iterations, it derives an optimal solution with a relatively small amount of additional time than for stochastic programs.

So, for the next large distribution system, we just adopt our customized bilinear Benders decomposition as the primary computational method. It is worth mentioning that
 the  substation capacity derived from the consideration of 100 scenarios is larger than those derived from other considerations, which definitely indicates that a reliable plan will be produced if sufficiently many stochastic scenarios should be respected. \\

\begin{table}[htbp]
  \renewcommand{\arraystretch}{1.0}
  \centering
  \caption{Load and price data}
  \label{table_I}
    \begin{tabular}{cccc}
    \toprule
    Time block & Load factor & Duration (hr) & Price (\$/MWh) \\
    \midrule
    1     & 0.65  & 2000  & 60 \\
    2     & 0.8   & 5760  & 70 \\
    3     & 1     & 1000  & 90 \\
    \bottomrule
    \end{tabular}%
  \label{tab:addlabel}%
\end{table}%

\begin{table}[htbp]
\renewcommand{\arraystretch}{1.1}
  \centering
  \caption{Stochastic scenario data}
  \label{table_II}
    \begin{tabular}{cccccc}
    \toprule
          & \multicolumn{5}{c}{Scenario scaling factor} \\
    \midrule
    scenario & 1     & 2     & 3     & 4     & 5 \\
    Case 1 & 2.02  & 0.97  & 2.8   & 2.89  & 0.76 \\
    Case 2 & 1.41  & 0.78  & 1.88  & 1.93  & 0.66 \\
        \midrule
    scenario & 6     & 7     & 8     & 9     & 10 \\

    Case 1 & 2.28  & 1.88  & 1.16  & 1.37  & 1.52 \\
    Case 2 & 1.57  & 1.33  & 0.89  & 1.02  & 1.11 \\
    \bottomrule
    \end{tabular}%
  \label{tab:addlabel}%
\end{table}%

\begin{table}[htbp]
\renewcommand{\arraystretch}{1.1}
  \centering
  \caption{Results of the chance-constrained model}
  \label{table_III}
    \begin{tabular}{cccc}
    \toprule
     \pbox{20cm} {Chance level \\1--$\epsilon$ (\%)} &  \pbox{20cm} {Discarded \\ scenarios} & Invested  branches &  \pbox{20cm} {Invested \\substation (pu)} \\
    \midrule
    100   & --      & 1--2, 2--5 & 0.11 \\
    90    & 4      & 1--2, 2--5 & 0.11 \\
    80    & 3, 4   & 1--2, 2--5 & 0.11 \\
    65    & 3, 4, 6 & 1--2, 2--5 & 0.11 \\
    \bottomrule
    \end{tabular}%
  \label{tab:addlabel}%
\end{table}
%
\begin{table}[htbp]
\renewcommand{\arraystretch}{1.1}
  \centering
  \caption{Computational performance of models}
  \label{table_IV}
    \begin{tabular}{ccccc}
    \toprule
    \pbox{20cm} {Scenario \\number} & \pbox{20cm} {Invested \\branch} & \pbox{20cm} {Substation\\ capacity (pu)} & \pbox{20cm} {Iteration \\number} & \pbox{20cm} {Time (sec.)\\ (gap \%)} \\
    \midrule
    \multicolumn{5}{c}{\textbf{Non-decomposed-SOCP}} \\
    10    & 1--2, 2--5 & 0.107 & --     & 43.0 \\
    50    & --     & --     & --     & T (N/A) \\
    100   & --     & --     & --     & T (N/A) \\
    \midrule
    \multicolumn{5}{c}{\textbf{Benders-linear StoP}} \\
    10    & 1--2, 2--5 & 0.11  & 3     & 23.2 \\
    50    & 1--2, 2--5 & 0.11  & 3     & 117.3 \\
    100   & 1--2, 2--5 & 0.118 & 3     & 241.0 \\
    \midrule
    \multicolumn{5}{c}{\textbf{Benders-SOCP StoP}} \\
    10    & 1--2, 2--5 & 0.11  & 4     & 5.9 \\
    50    & 1--2, 2--5 & 0.11  & 4     & 31.4 \\
    100   & 1--2, 2--5 & 0.118 & 4     & 72.2 \\
    \midrule
    \multicolumn{5}{c}{\textbf{Benders-SOCP CC (1--$\epsilon$ =90\%)}} \\
    10    & 1--2, 2--5 & 0.11  & 4     & 7.9 \\
    50    & 1--2, 2--5 & 0.11  & 4     & 50.9 \\
    100   & 1--2, 2--5 & 0.118 & 4     & 149.3 \\
    \bottomrule
    \end{tabular}%
  \label{tab:addlabel}\vspace{-10pt}%
\end{table}%

\vspace{-.5cm}
\subsection{18-node system}
\begin{figure}[t]
 \centering
\includegraphics[width=7cm, height=6cm]{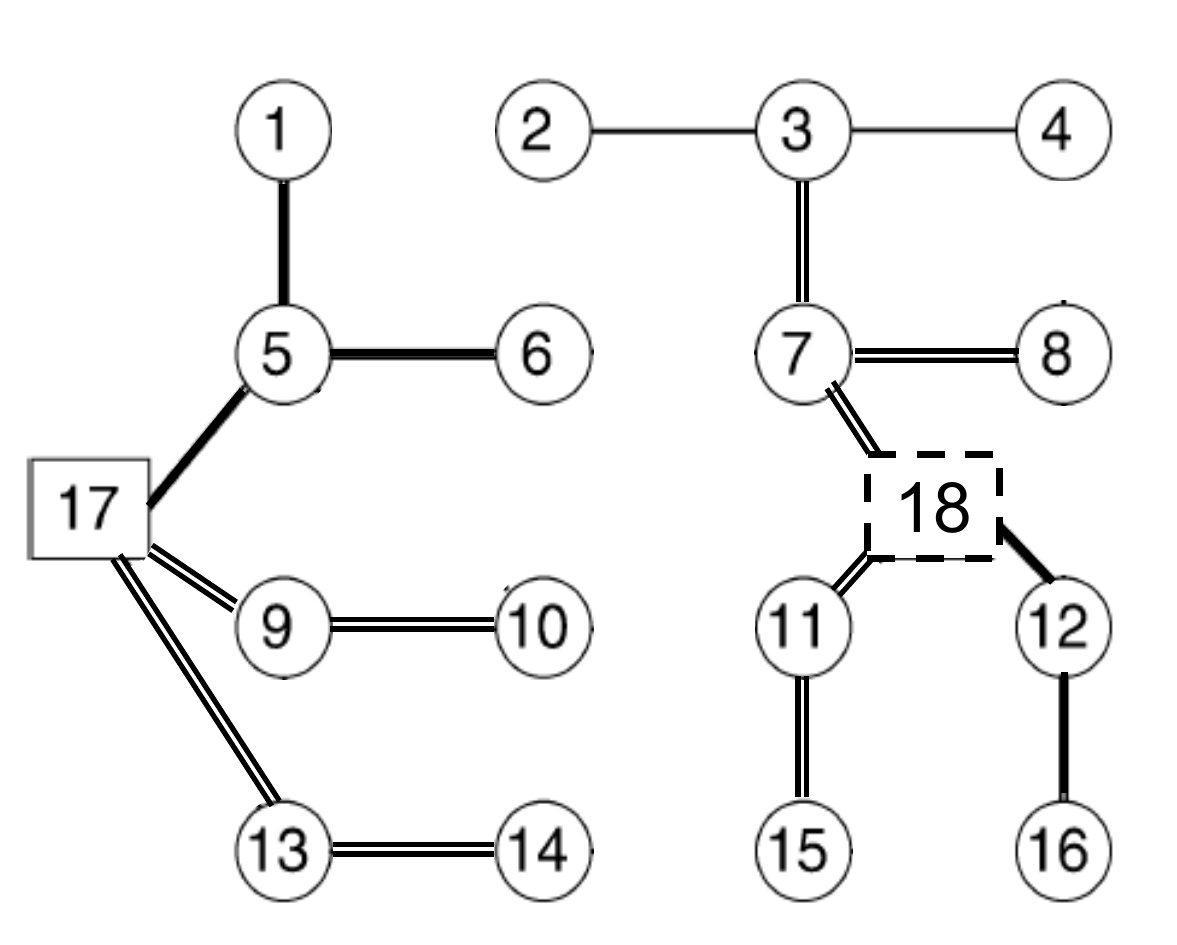}
\includegraphics[width=7cm, height=6cm]{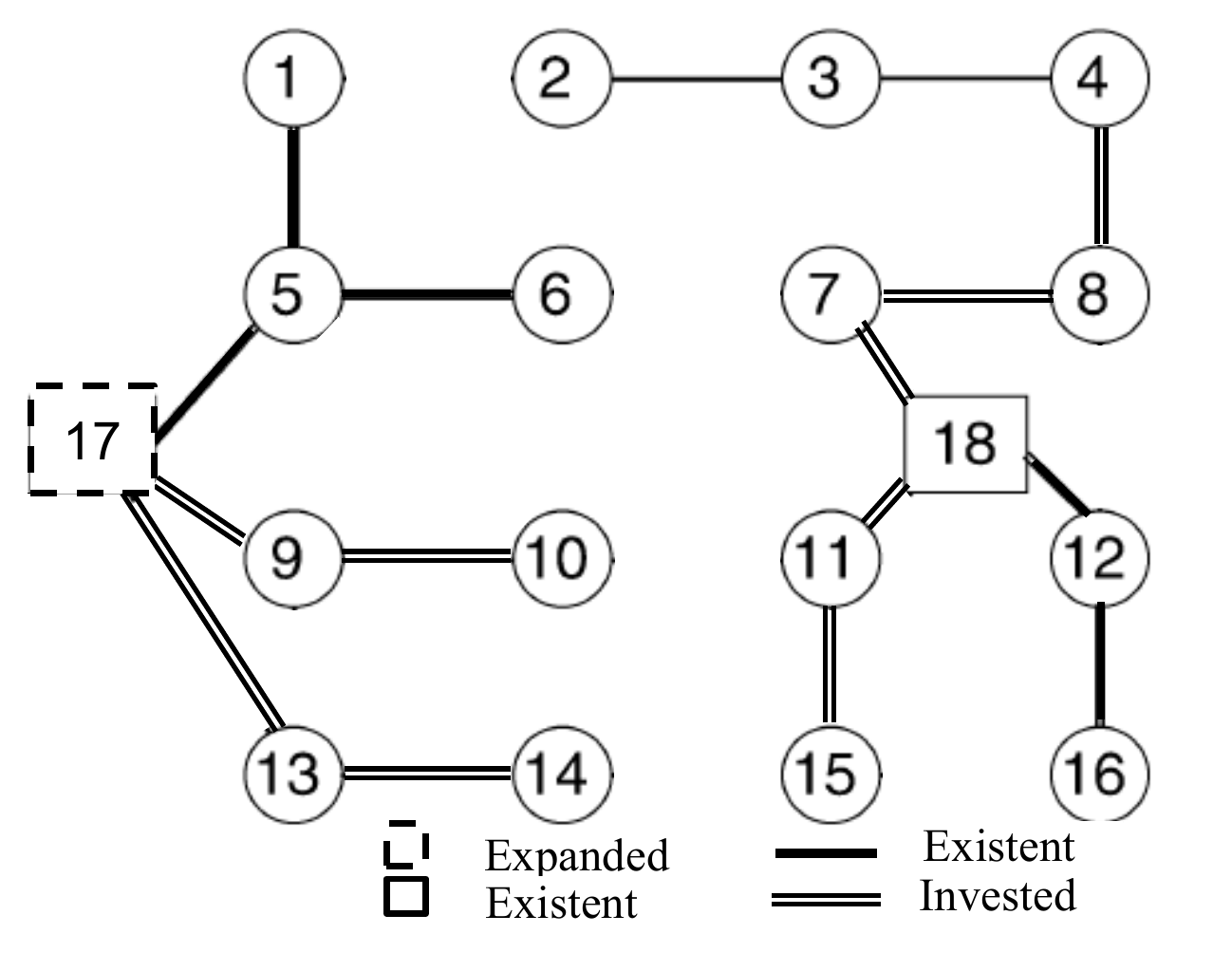}
\caption{Expansion results with chance of a) 100\% and b) 80\%}
\vspace{-10pt}
\label{fig_2.b}
\end{figure}

The second test system was adopted from \citep{Haffner2} with some modifications (see Fig.~\ref{fig_2.b}). This system has 18 nodes, 2 substations, and 24 branches. The existing and candidate branches are shown in Table~\ref{table_V}. The uncertainty of the price and load was modeled by twenty equiprobable scenarios with the scaling factor being uniformly distributed in the range [0.6, 1.8]. We considered three time blocks in the target year as shown Table~\ref{table_I} and assumed that every existing branch can be removed or re-wired.\par

\begin{table*}[htbp]
  \centering
  \caption{expansion results of the 18-node system}
  \label{table_V}
    \begin{tabular}{rrrrr}
    \toprule
    \multicolumn{1}{c}{\multirow{2}[2]{*}{Existent branch}} & \multicolumn{1}{c}{\multirow{2}[2]{*}{Candidate branch}} &
    \multicolumn{1}{c}{\multirow{2}[2]{*}{ \pbox{20cm} {Invested branch\\ (chance level 100\%)}}} & \multicolumn{1}{c}{\multirow{2}[2]{*}{\pbox{20cm} {Invested branch\\ (chance level 80\%)}}} &
    \multicolumn{1}{c}{\multirow{2}[2]{*}{\pbox{20cm} {Invested branch\\ (deterministic)}}}
    \\
    \multicolumn{1}{c}{} & \multicolumn{1}{c}{} & \multicolumn{1}{c}{} & \multicolumn{1}{c}{}\\
    \midrule
    1--2, 2--3 & 2--3, 4--8, 7--8, 9--13, 1--5, 5--6 & \multicolumn{1}{c}{7--8, 11--18} & \multicolumn{1}{c}{4--8, 7--8} & \multicolumn{1}{c}{4--8, 11--18}\\
     3--4, 1--5 & 10--11, 11--18, 13--14, 12--16  & \multicolumn{1}{c}{13--14, 13--17} & \multicolumn{1}{c}{11--18, 13--14} & \multicolumn{1}{c}{13--14, 13--17}\\
     5--6, 5--17  & 13--17, 14--15, 5--10, 3--7, 3--4 & \multicolumn{1}{c}{7--18, 9--10 } & \multicolumn{1}{c}{13--17, 7--18} & \multicolumn{1}{c}{7--18, 9--10}\\
        12--16, 12--18  & 1--2, 6--7, 7--18, 8--12, 15--16 & \multicolumn{1}{c}{9--17, 11--15 } & \multicolumn{1}{c}{9--10, 9--17} & \multicolumn{1}{c}{9--17, 11--15} \\
          & 9--10, 9--17, 11--15, 12--18, 5--17 & \multicolumn{1}{c}{3--7} & \multicolumn{1}{c}{11--15} & \multicolumn{1}{c}{3--7} \\
        \midrule
    \multicolumn{2}{r}{Substation investment (MVA)} & \multicolumn{1}{c}{2.7 (node 18)} & \multicolumn{1}{c}{2.7 (node 17)} & \multicolumn{1}{c}{--}\\
    \multicolumn{2}{r}{Total cost (\$)} & \multicolumn{1}{c}{444 136} & \multicolumn{1}{c}{387 331} & \multicolumn{1}{c}{205 784}\\
    \multicolumn{2}{r}{Time (sec.)} & \multicolumn{1}{c}{62.0} & \multicolumn{1}{c}{151.0} & \multicolumn{1}{c}{9.0}\\
    \bottomrule
    \end{tabular}
  \label{tab:addlabel}\vspace{-10pt}%
\end{table*}%

The problem was solved for target year 5 considering $(a)$ a deterministic case for a load level of 150\%; $(b)$ two stochastic cases with chance levels of 100 and 80 percent. Detailed results are provided in Table~\ref{table_V} and stochastic expansion results are depicted in Fig.~\ref{fig_2.b}. As seen from Table~\ref{table_V}, 9 new feeders are invested and added to the network and the substation is expanded by 2.7 MVA in the stochastic cases while the deterministic case does not invest in substation capacity. Clearly, different radial topologies are achieved with the costs given in Table~\ref{table_V} and the highest expansion cost is associated
with the case with 100 percent chance level, which is expected. Note that instances can be solved within 8-18 iterations, with solution times shown in the table.\\

\section {Conclusion}
In this paper, a stochastic second order mixed integer model was constructed and studied for distribution system expansion planning. A chance-constrained variant of the problem was also developed, which can be used to obtain a cost-effective plan by avoiding extreme scenarios. To solve the problem, we developed a novel bilinear Bender  decomposition  algorithm  that  handles the conic power flow formulations and deals a chance constraint imposed on stochastic scenarios.  On a set of  instances, we performed numerical experiments, analyzed our model's performance and generalized insights.

Based on our numerical experiments, it was demonstrated that our proposed algorithms have a remarkable strength, which drastically outperforms professional mixed integer conic solver CPLEX by orders of magnitude. They therefore greatly improves our solution capacity for this type of problems. Indeed,  to the best of our knowledge, it is the first time that chance-constrained stochastic mixed integer SOCP model can be solved efficiently. Because the developed algorithms are generic, we believe that they can be used in other application areas, such as power system expansion and operation planning.

\bibliographystyle{unsrt} 
\bibliography{zengbo} 
%

\end{document}